\documentclass[12pt,a4paper]{article}
\usepackage{a4wide}
\usepackage{amsmath,amsfonts,amssymb,amsthm}
\usepackage{mathrsfs}
\usepackage{ifthen}
\usepackage{epsfig}

\newtheorem{Theorem}{Theorem}[section]
\newtheorem{Cor}{Corollary}[section]
\newtheorem{Def}{Definition}[section]
\newtheorem{Lemma}{Lemma}[section]
\newtheorem{Prop}{Proposition}[section]

\makeatletter
\@addtoreset{equation}{section}
\makeatother

\makeatletter
\newcommand\txt[2]{{#1{\ifthenelse{\equal {#2}{}}{}{ #2}}. }}
\def\Remark#1.#2\par{\par\medskip\noindent\txt{\bf Remark}{#1}\ #2\par%
               \@ifnextchar {\Remark}{}{\@ifnextchar\begin{}{\smallskip}}}
\def\Example#1.{\txt{\bf Example}{#1}}
\makeatother

\newcommand\emp\varnothing
\newcommand\eps\varepsilon
\let\ds\displaystyle

\let\dd\partial

\def\dout{\dd^{\rm out}}

\def\^#1{^{\overline{#1}}}

\def\lchi{\mathop{\ell\chi}\nolimits}

\begin{document}

\title{Number of vertices in graphs with locally small chromatic number and large chromatic number}

\author{Ilya I. Bogdanov\thanks{The work was supported by the Russian government project 11.G34.31.0053 and by RFBR grant No.~13-01-00563.}}

\maketitle

\begin{abstract}
  We discuss the minimal number of vertices in a graph with a large chromatic number such that each ball of a fixed radius in it has a small chromatic number. It is shown that for every graph $G$ on $\sim((n+rc)/(c+rc))^{r+1}$ vertices such that each ball of radius $r$ is properly $c$-colorable, we have $\chi(G)\leq n$.
\end{abstract}

\section{Introduction}

Let $G=(V,E)$ be a graph (with no loops or multiple edges). By $d(u,v)=d_G(u,v)$ we denote the {\em distance} between the vertices $u,v\in V$. A subset $V_1\subseteq V$ is {\em independent} if none of the edges has both endpoints in $V_1$. The {\em chromatic number} $\chi(G)$ of~$G$ is the minimal number of colors in a proper coloring of~$G$, that is --- the minimal number of parts in a partition of $V$ into independent subsets.

\begin{Def}
  Let $r$ be a nonnegative integer. The {\em ball} of radius $r$ with center $v\in V$ is the set $U_r(v,G)=\{u\in G\colon d(u,v)\leq r\}$. For $r\geq 1$, the {\em $r$-local chromatic number} $\lchi_r(G)$ of a graph~$G$ is the maximal chromatic number of a ball of radius~$r$ in~$G$.
\end{Def}

Notice that even for $r=1$ our definition of the local chromatic number is quite different from that introduced by Erd\H os et al. in~\cite{erd-loc}.

By a well-known result of Erd\H os~\cite{erdos}, for every integer $g>2$ and every $n>n_0(g)$, there exists a graph on $n^{2g+1}$ vertices of girth~$g$ and chromatic number greater than~$n$; thus for every $r$ there exist a graph~$G$ with $\lchi_r(G)=2$ and arbitrarily large $\chi(G)$. Later Erd\H os~\cite{erdos3} conjectured that for every positive integer~$s$ there exists a constant $c_s$ such that the chromatic number of each graph~$G$ having $N$ vertices and containing no odd cycles of length less than $c_sN^{1/s}$ does not exceed $s+1$. This conjecture was proved by Kierstead, Szemer\'edi, and Trotter~\cite{kst}. In fact, they have proved a more general result; we  will formulate this result in terms of the following notion.

\begin{Def}
  Let $n$, $r$, and $c$ be positive integers. Denote by~$f_c(n,r)$ the maximal integer~$f$ with the following property: If $G$ is a graph on $f$ vertices and $\lchi_r(G)\leq c$ then $\chi(G)\leq n$.
\end{Def}

Then the aforementioned result can be formulated as
\begin{equation}
  f_c\bigl(k(c-1)+1,r\bigr)\geq\left\lfloor\frac r{2k}\right\rfloor^k,
  \label{lower-kst}
\end{equation}
while the result by Erd\H os yields
\begin{equation}
  f_c(n,r)\leq f_2(n,r) <n^{4r+5} \qquad \text{for all $n>n_0(r)$.} 
  \label{upper-erd}
\end{equation}

Several examples (cf., for instance, \cite{schrijver,stie}) show that the estimate~\eqref{lower-kst} has a correct order in $r$  for $c=2$. In~\cite{bogd-nrv} we show that this order is sharp even for all $c$. namely, it is shown that
$$
  f_c\bigl(k(c-1),r\bigr)<\frac{(2rc+1)^k-1}{2r}.
$$

On the other hand, the estimate~\eqref{lower-kst} does not work for $n\gtrsim (c-1)r$. For $c=2$, Berlov and the author~\cite{berl-bogd} obtained the estimate
\begin{equation}
  \label{bb}
  f_2(n,r)\geq \frac{(n+r+1)(n+r+2)\cdots(n+2r+1)}{2^r(r+1)^{r+1}}.
\end{equation}

It is worth mentioning that for several specific series of parameters there exist almost tight bounds of $f_c(n,r)$. Firstly, asymptotics of $f_2(n,1)$ is tightly connected with the asymptotics of Ramsey numbers $R(n,3)$. In the papers of Ajtai, Koml\'os, and Szemer\'edi~\cite{aks} and Kim~\cite{kim} it is shown that $c_1\frac{n^2}{\log n}\leq
R(n,3)\leq c_2\frac{n^2}{\log n}$ for some absolute constants $c_1,c_2$. One can check that these results imply the bounds
$$
  c_3 n^2\log n\leq f(n,2)\leq c_4 n^2 \log n
$$
for some absolute constants $c_3,c_4$.

The asymptotics of $f_2(3,r)$ is also well investigated. From the generalized Mycielski construction by Stiebitz~\cite{stie} it follows that $f_2(3,r)<2r^2+5r+4$. On the other hand, Jiang~~\cite{jiang} showed that $f_2(3,r)\geq (r-1)^2$.

The aim of this paper is to extend this estimate for larger values of $c$. For the convenience, we use the notation $n\^k=n(n+1)\dots(n+k-1)$. We prove the following results.

\begin{Theorem}
  \label{th-gen}
  For all positive integer $n$, $r$, and $c>1$ we have
  \begin{equation}
    f_c(n,r)\geq \frac{(n/c+r/2)\^{r+1}}{(r+1)^{r+1}}.
    \label{eq-gen}
  \end{equation}
\end{Theorem}

%%\begin{Theorem}
%%  \label{th-large}
%%\end{Theorem}

\section{Main result}
\label{sec-low}

For a graph $G=(V,E)$ and a subset $V_1\subseteq V$, we denote by $G[V_1]$ the induced subgraph of $G$ on the set $V_1$. For $v\in V$, we denote by $S_r(v,G)=\{u\in V\mid d_G(u,v)= r\}$ the {\em sphere} with radius~$r$ and center~$v$. In particular, $S_0(v,G)=U_0(v,G)=\{v\}$. Denote also by $\dout_G V_1=\{u\in V\setminus V_1\mid \exists v\in V_1: (u,v)\in E\}$ the {\em outer boundary} of a subset $V_1\subseteq V$. In particular, $S_r(v,G)=\dout_G U_{r-1}(v,G)$.

Our estimate is based on the following lemma.

\begin{Lemma}
  For every graph $G=(V,E)$ and every positive integer $r$, there exists a decomposition $V=U\sqcup N$ such that each connected component of $U$ lies in some ball in $G$ of radius $r$, and 
  \begin{equation}
    |N|\leq \frac{\sqrt[r+1]{|V|}-1}{\sqrt[r+1]{|V|}}|V|.
    \label{cond-ind}
  \end{equation}
  \label{lm-main}
\end{Lemma}

\proof
Set $v=|V|$. We will construct inductively a sequence of partitions of~$V$ into nonintersecting parts,
$$
  V=U_1\sqcup U_2\sqcup\dots\sqcup U_s\sqcup N_s\sqcup V_s,
$$
such that the following conditions are satisfied:

(i) for all $i=1,\dots,s$ we have $\dout_G U_i\subseteq N_s$; moreover, $\dout_G V_s\subseteq N_s$;

(ii) for every $i=1,2,\dots,s$ the graph $G[U_i]$ is contained in some ball in $G$ of radius $r$;

(iii) $(\sqrt[r+1]{v}-1)(|U_1|+\dots+|U_s|)\geq |N_s|$.

For the base case $s=0$, we may set $V_0=V$, $N_0=\emp$ (there are no sets~$U_i$ in this case).

For the induction step, suppose that the partition $V=U_1\sqcup U_2\sqcup\dots\sqcup U_{s-1}\sqcup N_{s-1}\sqcup V_{s-1}$ has been constructed, and assume that the set~$V_{s-1}$ is nonempty. Consider the graph $G_{s-1}=G[V_{s-1}]$ and choose an arbitrary vertex $u\in V_{s-1}$. Now consider the sets
$$
  U_0(u,G_{s-1})=\{u\},\quad  U_1(u,G_{s-1}), \quad \dots, \quad U_{r+1}(u,G_{s-1}).
$$
One of the ratios
$$
  \frac{|U_1(u,G_{s-1})|}{|U_0(u,G_{s-1})|}, \quad 
  \frac{|U_2(u,G_{s-1})|}{|U_1(u,G_{s-1})|}, \quad
  \dots, \quad 
  \frac{|U_{r+1}(u,G_{s-1})|}{|U_r(u,G_{s-1})|}
$$
does not exceed~$\sqrt[r+1]v$, since the product of these ratios is
$$
  |U_{r+1}(u,G_{s-1})|\leq v.
$$
So, let us choose $1\leq m\leq r+1$ such that
$$
  \frac{|U_m(u,G_{s-1})|}{|U_{m-1}(u,G_{s-1})|}\leq \sqrt[r+1]v.
$$

Now we set
$$
  U_s=U_{m-1}(u,G_{s-1}), \quad N_s=N_{s-1}\cup S_m(u,G_{s-1}),
  \quad
  V_s=V_{s-1}\setminus U_m(u,G_{s-1}).
$$
Since the condition~(i) was satisfied on the previous step, we have
$$
  (\dout_G V_s)\cup(\dout_G U_s)\subseteq \dout_G V_{s-1}\cup S_m(u,G_{s-1})\subseteq N_s
$$
so this condition also holds now. The condition~(ii) is satisfied trivially. Finally, the choice of~$m$ and the condition~(iii) for the previous step imply that
\begin{gather*}
  \sqrt[r+1]{v}\cdot |U_s|=\sqrt[r+1]{v}\cdot |U_{m-1}(u,G_{s-1})|\geq |U_m(u,G_{s-1})|,\\
  (\sqrt[r+1]{v}-1)(|U_1|+\dots+|U_{s-1}|)\geq |N_{s-1}|
\end{gather*}
and hence
$$
  (\sqrt[r+1]{v}-1)(|U_1|+\dots+|U_{s}|)\geq
  |N_{s-1}|+|U_m(u,G_{s-1})|-|U_{m-1}(u,G_{s-1})|
  =|N_s|.
$$
Thus, the condition~(iii) also holds on this step.

Continuing the construction in this manner, we will eventually come to the partition with $V_s=\emp$ since the value of $|V_s|$ strictly decreases. As the result, we obtain the partition $V=U_1\sqcup U_2\sqcup\dots\sqcup U_s\sqcup N_s$ such that $|N_s|\leq
(\sqrt[r+1]{v}-1)(|U_1|+\dots+|U_s|)$. So, setting $U=U_1\cup\dots\cup U_s$ and $N=N_s$ we get
$$
  \sqrt[r+1]{v}\cdot |N|\leq (\sqrt[r+1]{v}-1)|U|+(\sqrt[r+1]{v}-1)|N|=|V|(\sqrt[r+1]{v}-1),
$$
or $\ds |N|\leq |V|\frac{\sqrt[r+1]{v}-1}{\sqrt[r+1]{v}}$, as required.
\qed

\begin{Cor}
  Setting $v=f_c(n,r)+1$, we have
  \begin{equation}
    v\geq \frac{\sqrt[r+1]{v}}{\sqrt[r+1]{v}-1}\bigl(f_c(n-2,k)+1\bigr).
    \label{eq-cor}
  \end{equation}
  \label{cor-est}
\end{Cor}

\proof Let $G=(V,E)$ be a graph on $v$ vertices such that $\lchi_r(G)\leq c$ although $\chi(G)>n$. Applying Lemma~\ref{lm-main} we get a decomposition $V=U\sqcup N$ such that $G[U]$ has a proper coloring in $c$ colors. So, $G[N]$ cannot be properly colored in $n-c$ colors, hence $|N|\geq f_c(n-c,r)$, hence the relation~\eqref{cond-ind} yields~\eqref{eq-cor}.
\qed

The next proposition shows how to make an explicit estimate for $f_c(n,r)$ from Corollary~\ref{cor-est}. 

\begin{Prop}
  Suppose that for some integer $n_0\geq 1$ and real $a$ the inequality
  \begin{equation}
    f(m,k)\geq \frac{(a+m/c)\^{r+1}}{(r+1)^{r+1}}-1
    \label{asym1}
  \end{equation}
  holds for $m=n_0$. Then the same estimate holds for all integer $m\geq n_0$ with $m-n_0\equiv 0\pmod c$.
  \label{recurr}
\end{Prop}

\proof
We use the Induction on $m$ with step $c$; the base case holds by the conditions of the proposition. 

Assume now that \eqref{asym1} holds for $m=n-c$ but not for $m=n$. Denote $v=f_c(n,r)+1$; then we have
$$
  \sqrt[r+1]v<\sqrt[r+1]{\frac{(a+n/c+r)^{r+1}}{(r+1)^{r+1}}}=\frac{a+n/c+r}{r+1}
$$
and hence
$$
  \frac{\sqrt[r+1]v}{\sqrt[r+1]v-1}>\frac{a+n/c+r}{a+n/c-1}.
$$
By~\eqref{eq-cor} this yields
$$
  v\geq \frac{a+n/c+r}{a+n/c-1}\cdot \frac{(a+(n-c)/c)\^{r+1}}{(r+1)^{r+1}}
  =\frac{(a+n/c)\^{r+1}}{(r+1)^{r+1}},
$$
which contradicts our assumption. Thus the induction step is proved.
\qed

\medskip
\textit{Proof of Theorem~\ref{th-gen}.}
In view of Proposition~\ref{recurr} it suffices to check~\eqref{eq-gen} for all $n\leq c$. Trivially, we have $f_c(n,r)\geq 1$. On the other hand, by the AM--GM inequality we have
$$
  \left(\frac nc+\frac r2\right)\^{r+1}\leq\left(\frac nc+r\right)^{r+1}\leq (r+1)^{r+1},
$$
thus 
$$
  f_c(n,r)\geq 1\geq \frac{(n/c+r/2)\^{r+1}}{(r+1)^{r+1}},
$$
as required.
\qed

\Remark. For large values of parameters, one may use the larger value of $a$ in Proposition~\ref{recurr}, for instance, by the use of~\eqref{lower-kst}.


\begin{thebibliography}{99}

\bibitem{aks} M. Ajtai, J. Koml\'os, and E. Szemer\'edi, A note on Ramsey numbers, J. Combin. Theory A, 29 (1980), 354--360.

\bibitem{bogd-nrv} I.I. Bogdanov. Examples of topologically highly chromatic graphs with locally small chromatic number. // arXiv:1311.2844
%
\bibitem{erdos} P. Erd\H os, Graph theory and probability. // Canad. J. Math., 11 (1959), 34--38.

\bibitem{erdos3} P. Erd\H os, Problems and results in graph theory and
combinatorial analysis. // In: Graph Theory and Related Topics, Academic Press,
New York, 1979, 153--163.

\bibitem{erd-loc} P. Erd\H os, Z. F\"uredi, A. Hajnal, P. Komj\'ath, V. R\"odl, \'A. Seress, Coloring graphs with locally few colors. // Discrete Math., 59 (1986), 21--34.

\bibitem{kst} H.~A.~Kierstead, E. Szemer\'edi, and W.~T.~Trotter, On coloring graphs with locally small chromatic number. // Combinatorica 4 (1984), 183--185.

\bibitem{jiang} T. Jiang. Small odd cycles in 4-chromatic graphs. // J. Graph Theory, 37 (2001), 115--117.

\bibitem{kim} J.~H.~Kim, The Ramsey Number $R(3,t)$ has order of magnitude $t^2/\log t$,
Random Structures and Algorithms, 7 (1995), 173–-207.

\bibitem{lovasz} L. Lov\'asz, Kneser's conjecture, chromatic number, and homotopy. // J. of Comb. Theory A 25 (1978), 319--324.

%\bibitem{SS} H. Sachs, M. Stiebitz, On constructive methodsin the theory of colour-critical graphs, Discrete Math. 74 (1989), 201--226.

\bibitem{schrijver} A. Schrijver, Vertex-critical subgraphs of Kneser graphs. // Nieuw Archief voor Wiskunde 26 (1978), 454--461.

\bibitem{stie} M. Stiebitz, Beitr\"age zur Theorie der f\"arbungskritischen Graphen. Habilitation, TH Ilmenau, 1985.

\bibitem{berl-bogd} S.~L.~Berlov and I.~I.~Bogdanov, On graphs with a large chromatic number that contain no small odd cycles. // J. Math. Sci. 184(2012), 573--578.

\end{thebibliography}
\end{document}